\documentclass[a4paper,12pt]{article}
\usepackage{amssymb,latexsym,amsmath,enumerate,verbatim,amsfonts,amsthm}

\def\tto{\;{\lower 1pt\hbox{$\rightarrow$}}\kern -10pt
\hbox{\raise 2pt\hbox{$\rightarrow$}}\;}

\def\ra{\rangle}
\def\la{\langle}

\def\R{\mathbb R}

\def\co{\mbox{\rm co}\,}

\def\dom{\mbox{\rm dom}\,}

\def\cl*co{\mbox{\rm cl}^*\mbox{\rm co}\,}
\def\cl{\mbox{\rm cl}\,}

\def\hs7{\hspace*{7pt}}

\newtheorem{thm}{Theorem}[section]

\newtheorem{rem}{Remark}[section]
\newtheorem{defi}{Definition}[section]
\newtheorem{cor}{Corollary}[section]
\newtheorem{ex}{Example}[section]

\newtheorem{lem}{Lemma}[section]

\newcommand{\bc}{\begin{center}}
\newcommand{\ec}{\end{center}}

\newcommand{\lrar}{\longrightarrow}

\newcommand{\be}{\begin{equation}}
\newcommand{\ee}{\end{equation}}
\newcommand{\beqa}{\begin{eqnarray*}}
\newcommand{\eeqa}{\end{eqnarray*}}

\newcommand{\mes}{\medskip}

\newcommand \pr {\noindent{\it Proof: \/\ }}
\newcommand{\lan}{\langle}
\newcommand{\ran}{\rangle}

 \newcommand{\p}{\partial}

\newcommand{\tK}{\tilde{K}}

 \newcommand{\gl}{\lambda}

\newcommand{\ga}{\alpha}

\newcommand{\noin}{\noindent}

\def\cmark{\mbox{$\rm\bf\rule{0.06em}{1.45ex}\kern-0.3em C$}}
\def\nmark{\mbox{$\rm\bf\rule{0.06em}{1.45ex}\kern-0.05em N$}}
\def\rmark{\mbox{$\rm\bf\rule{0.06em}{1.45ex}\kern-0.05em R$}}
\normalsize

\begin{document}

\title{{{{{\textsf{Constrained Best Approximation with  Nonsmooth Nonconvex Constraints}}}}}}
\date{}
\author{H. Mohebi    \thanks{Department of Mathematics and Mahani Mathematical Research Center,
 Shahid Bahonar University of Kerman, P.O. Box: 76169133, Postal Code: 7616914111, Kerman, Iran, Phone $\&$ Fax: +98-34-33257280,   e-mail: hmohebi@uk.ac.ir  (Hossein Mohebi)}  \thanks{This  work  was partially  performed while  the author was on sabbatical leave  at the Department of Applied Mathematics, University of New South Wales, Sydney,  NSW, 2052, Australia.}}
\maketitle

\mes

\begin{abstract} \noin In this paper, we consider the constraint set $K$ of inequalities with  nonsmooth nonconvex constraint functions.  We show  that under Abadie's constraint qualification the   "perturbation property" of the best approximation to any $x$ in $\R^n$ from a convex set $\tK:=C \cap K$ is characterized by the strong conical hull intersection property (strong CHIP) of $C$ and $K,$  where $C$ is a non-empty closed convex subset of $\R^n$  and the set $K$ is represented by $K:=\{x\in \R^n : g_j(x) \le 0, \ \forall \ j=1,2,\ldots,m \}$  with  $g_j : \R^n \lrar \R$ $(j=1,2, \cdots,m)$  is a tangentially convex function at a given  point $\bar x \in K.$  By using the idea
of tangential subdifferential  and a non-smooth version of  Abadie's constraint qualification, we do this by first proving    a dual cone characterization of the constraint set $K.$  Moreover, we present sufficient conditions for which the strong CHIP property holds.  In particular, when the set $\tK$ is  closed and convex, we show that the Lagrange multiplier characterization of best approximation holds under a non-smooth version of   Abadie's constraint qualification. The obtained results extend many  corresponding results in the context of constrained best approximation.  Several examples are provided to clarify the results.

  \vskip 0.1cm%

\vskip 0.5cm\vskip 0.1cm\noin \textbf{2010 (AMS) Mathematics Subject Classification.}  41A29; 41A50; 90C26; 90C46.
\end{abstract}

\mes

\mes

\section{Introduction}
\noin  Finding suitable conditions for "perturbation property" has been of substantial interest in constrained best approximation as it is often easier to compute the best approximation from a closed convex set $C$ than from the constraint set $K.$  The merit and  motivation  for such characterization $($perturbation property$)$ is inspired from \cite[Chapter 10]{DB}.  "Characterizing constrained interpolation from a convex set" is one of the applications of  the "perturbation property" $($for more details, see \cite{DB}$).$

\mes
\noin For many years, a great deal of attention has been focusing on the case where the constraint set $K:=\{x \in \R^n : g_j(x) \le 0, \ \forall \ j=1,2,\ldots,m \}$ is a closed convex set and has a convex representation in the sense that  $g_j$ $(j=1,2,\cdots,m)$  is a convex function \cite{DB,DLW1,DLW2,jm1,LJ,LN1}. Various characterizations of the perturbation property
have been given by using local constraint qualifications such as
the strong conical hull intersection property (strong CHIP) of $C$ and $K$ at the best approximation \cite{DE,DLW2,jm1,LJ}.

\mes

\noin In this paper, we study the problem of whether the best approximation to any $x \in \R^n$
from the {\it closed convex} set $\tK:=C \cap K$  can be characterized by the best approximation to a perturbation $x- x^*$ of $x$ from  a closed convex set $C \subseteq \R^n$ for some $x^*$ in a certain cone in $\R^n,$  where $K:=\{x \in \R^n : g_j(x) \le 0, \ \forall \ j=1,2,\ldots,m \}$ with  $g_j : \R^n \lrar \R$ $(j=1,2,\cdots,m)$ is  a tangentially convex function at the best approximation.  We show that the strong CHIP of $C$ and $K$ at the best approximation continues to completely characterize the perturbation property of the best approximation from the closed convex set $\tK:=C \cap K$ under  Abadie's constraint qualification.

\mes
\noin Indeed, by using the idea of tangential subdifferential  and a  non-smooth version of  Abadie's constraint qualification (which is the weakest qualification among the other well known constraint qualifications), we prove this by first establishing a dual cone characterization of the constraint set $K.$  In  the special case, when $\tK$ is a  closed convex set, we show that the Lagrange multipliers characterization of best approximation holds under a non-smooth version of  Abadie's constraint qualification.   Our results  recapture the corresponding known results of \cite{DE,DLW1,DLW2,jm1,jm2,jm3,LJ,LN1,LN2,ms}. Several illustrative examples are presented to clarify our results.

\mes
\noin The paper has the following structure. In Section 2, we provide the basic results on tangentially convex functions and a non-smooth version of the constraint qualifications.  Dual cone characterizations of the constraint set $K$  and sufficient conditions for which the strong CHIP holds are presented in Section 3. In Section 4, we first show that the strong CHIP  completely characterizes the perturbation property of the best approximation. Finally, we show that under a non-smooth version of  Abadie's constraint qualification the Lagrange multipliers characterization of best approximation holds. Also, several  examples are presented to illustrate  our results.

\mes

\mes

\section{Preliminaries}
\noin We start this section by fixing notations and preliminaries that will be used later. Recall \cite{bor} that  for a function $f:\R^n \lrar \R,$ the directional derivative of $f$ at a point $\bar x \in \R^n$  in the direction $\nu \in \R^n$ is defined by   \be \label{41} f'(\bar x, \nu) :=\lim_{\ga \lrar 0^+} \frac{f(\bar x + \ga \nu) - f(\bar x)}{\ga}, \ee if the limit exists.
\noin Recall \cite{ml1} that  a function $f:\R^n \lrar \R$  is called tangentially convex at a point $\bar x \in \R^n,$  if  $f'(\bar x, \cdot)$ is  a real valued  convex function.

\mes
\noin It should be noted that if the function $f$ is tangentially convex at a point $\bar x \in\R^n,$ then, since $f'(\bar x, \cdot)$ is a positively homogeneous function, we conclude that $f'(\bar x, \cdot)$ is  a  sublinear function on $\R^n.$

\mes

\noin The tangential subdifferential of a function $f:\R^n \lrar \R$ at a point $\bar x \in \R^n$ is defined by \be \label{42} \p_T f(\bar x) :=\{x^* \in \R^n : \lan x^*, \nu \ran \le f'(\bar x, \nu), \ \forall \ \nu \in \R^n\}. \ee

\noin If $f$ is tangentially convex at $\bar x,$ then, $\p_T f(\bar x) \neq \emptyset,$ and moreover, $f'(\bar x, \cdot)$ is the support functional of $\p_T f(\bar x),$ i.e., for each $\nu \in \R^n,$ we have  \be \label{4.35} f'(\bar x, \nu) = \max_{x^* \in \p_T f(\bar x)} \lan x^*, \nu \ran. \ee
\noin It should be noted that if $f$ is a convex function, then, $\p_T f(x) = \p f(x)$ for each $x \in \R^n,$ where $\p f(x)$ is the classical convex subdifferential of $f$ at $x.$

\mes

\begin{rem} \label{rem2.1} Note that if the function $f:\R^n \lrar \R$  is  tangentially convex at a point $\bar x \in \R^n,$  then,  $f'(\bar x, \cdot)$ is  a real valued  convex function   on $\R^n,$  and hence, $f'(\bar x, \cdot)$ is a continuous function on $\R^n.$ \end{rem}

\mes

\noin Now, let $K \subseteq \R^n$  be defined by
\be \label{4.1} K:= \{x \in \R^n:  g_j(x) \le 0, \ \forall \ j=1,2,\ldots,m \},\ee where  $g_j : \R^n \lrar \R$ $(j=1,2, \cdots,m)$  is a tangentially convex function at a given  point $\bar x \in K.$  Let $C$ be a non-empty closed convex subset of $\R^n$ such that $C \cap K \neq \emptyset,$ and let  $S:=\R_+^m.$ Note that $K$ is   not necessarily a closed or a convex set. Let  \be \label{5000} \tK := C \cap K, \ee    and \be \label{4.00} I :=\{1,2,\cdots,m\}. \ee   \noin  For a point $\bar x \in K,$   we define \be \label{4.000} I(\bar x) :=\{j \in  I : g_j(\bar x) = 0 \}. \ee

\mes
\noin For a set $W \subseteq \R^n,$  let \be \label{0q2}  W^\circ:=\{\lambda\in\R^n\,:\, \la \lambda, y\ra\leq 0, \ \forall \  y\in W\}, \ee  where we denote $\lan \cdot, \cdot \ran$  for the inner product of $\R^n.$
\noin The normal cone  to a set $H \subseteq \R^n$ at a point $x \in \R^n$ is defined by \be \label{1.1}  N_H(x) := \{u \in \R^n :  \lan u, t - x \ran \le 0, \ \forall \ t \in H \}. \ee  It is clear that \beqa N_H(x) = (H - x)^\circ, \  \ (x \in \R^n). \eeqa
\noin  Let $U$ be a subset of $\R^n,$ and  let $x \in \R^n.$ We recall \cite{BS,bor} the contingent cone  of $U$ at $x$ is defined by \begin{eqnarray} \label{1.2} &&T_U(x) \nonumber \\ &:=&\{x^*\in \R^n : \exists \ \{\ga_k\}_{k \ge 1} \subset \R_{++} \ \mbox{and} \ \exists \ \{x_k^*\}_{k \ge 1} \subset \R^n \ \mbox{such that} \nonumber \\ &&\ga_k \lrar 0^+, \  x_k^* \lrar x^* \ \mbox{and} \ x + \ga_k x_k^* \in U, \ \forall \ k \ge 1 \}. \end{eqnarray}
\noin  We now introduce a non-smooth version of the linearized tangential cone: \be \label{1.3}  D(\bar x):=\{x^* \in \R^n : \lan x^*, \eta_j \ran \le 0, \ \forall \ \eta_j \in \p_T g_j(\bar x), \ \forall \ j \in I(\bar x) \}, \ee where $\bar x \in K.$

\mes
\noin Note that the  non-smooth  linearized tangential cone $D(\bar x)$  reduces to  its counterpart in  the case of differentiability \cite{BS,bor}. Moreover, $D(\bar x)$ is a convex cone.

\mes
\noin We now present the definition of  the nearly convexity which has been given in \cite{jm3}. Let $V$ be a non-empty subset of $\R^n$ and $x \in V.$

\mes
\begin{defi} \label{def00}  {\bf $($Nearly Convex at $x\in V).$} The set  $V$  is nearly convex at the point  $x \in V$ if for each  $y \in V$ there exists a  sequence   $\{t_k\}_{k\ge 1}$  of positive real numbers   with  $t_k \lrar 0^+$  such that  $x+ t_k(y-x) \in V$ for all sufficiently large  $k \in \mathbb{N}.$  \\
\noin The set $V$ is called {\it nearly convex} whenever it is nearly convex at each of its points. It is easy to check that if $V$ is convex, then it is nearly convex at each $x\in V.$   As shown in \cite{h1}, the nearly convexity may hold at a point for a non-convex set $($for more details and illustrative examples related to  the nearly convexity, see \cite{h1,jm3}$).$  \end{defi}

\mes
\begin{lem} \label{lem10}  Let $K$ be closed, given by $(\ref{4.1}),$ and let $C$ be a non-empty closed convex subset of $\R^n$ such that $C \cap K \neq \emptyset.$ Let $\tK := C \cap K,$ and   $\bar x\in \tK.$   Assume that $K$ is nearly convex at the point $\bar x.$   Then,  $T_{\tK}(\bar x) \subseteq D(\bar x),$ where $T_{\tK}(\bar x)$ and  $D(\bar x)$ defined by $(\ref{1.2})$ and $(\ref{1.3}),$ respectively.   \end{lem} \pr Let $x^* \in T_{\tK}(\bar x)$ be arbitrary. Then there exist sequences $\{\ga_k\}_{k \ge 1} \subset \R_{++}$ and $\{x_k^*\}_{k \ge 1} \subset \R^n$  such that $\ga_k \lrar 0^+,$  $x_k^* \lrar x^*$  and  $\bar x + \ga_k x_k^* \in \tK$ for all $k \ge 1.$  Since, by the hypothesis, $K$  is nearly convex at the point $\bar x$  and $\bar x + \ga_k x_k^* \in K$ for all $k \ge 1,$    it follows from  Definition \ref{def00}  that, for each $k \ge 1,$   there exists a sequence $\{\beta_{k, p}\}_{p \ge 1} \subset \R_{++}$ with $\beta_{k, p} \lrar 0^{+}$ $($as $p \lrar +\infty)$  such that  $\bar x + \beta_{k, p}(\bar x + \ga_k x_k^*- \bar x) \in K$  for all  sufficiently large $p \in \mathbb{N}.$ This implies that \be \label{1.4}  g_j(\bar x + \beta_{k, p} \ga_k x_k^*) \le 0, \ \mbox{for all sufficiently large} \ p \in \mathbb{N}, \ \forall \ k \ge 1, \  \forall \ j \in I. \ee  Since $g_j$ $(j \in I)$ is tangentially convex at $\bar x,$  it follows, by the definition,  that $g_j^\prime(\bar x, \cdot)$ is a real valued positively homogeneous and convex function on $\R^n.$    Therefore, for each $j \in I(\bar x),$  in view of $(\ref{1.4})$ we have \beqa g_j^\prime(\bar x,  \ga_k x_k^*) &=& \lim_{p \lrar +\infty} \frac{g_j(\bar x + \beta_{k, p} \ga_k x_k^*) - g_j(\bar x)}{\beta_{k, p}} \\ &=& \lim_{p \lrar +\infty}  \frac{g_j(\bar x + \beta_{k, p} \ga_k x_k^*)}{\beta_{k, p}}  \\ &\le&0, \ \forall \ k \ge 1, \eeqa  and hence,  \be \label{1.5}
 g_j^\prime(\bar x,   x_k^*) \le 0, \ \forall \ k \ge 1. \ee
Since $x_k^* \lrar x^*$  and $g_j^\prime(\bar x, \cdot)$ is continuous on $\R^n$ $($see Remark \ref{rem2.1}$),$   we conclude from (\ref{1.5}) that \beqa g_j^\prime(\bar x, x^*) \le 0, \ \forall \ j \in I(\bar x). \eeqa  This together with $(\ref{4.35})$ implies that \beqa \lan x^*, \eta_j \ran \le 0, \ \forall \ \eta_j \in \p_T g_j(\bar x), \ \forall \ j \in I(\bar x), \eeqa and so, $x^* \in D(\bar x),$ which completes the proof.  \hfill \rule{2mm}{2mm}

\mes
\noin Now, let us  to define the non-smooth versions of
Robinson's constraint qualification and  Abadie's constraint qualification.

\mes

\begin{defi} \label{def2.2} {\bf $($Non-smooth Version of   Robinson's Constraint Qualification $(NRCQ)).$} Let  $K =\{x \in \R^n : g_j(x) \le 0, \ \forall \ j=1,2,\cdots, m\}$ be  as in $($\ref{4.1}$),$  and let $C$ be a non-empty closed convex subset of $\R^n$ such that $C \cap K \neq \emptyset.$ Let $\tK := C \cap K,$ and   $\bar x\in \tK.$  We say that  non-smooth   Robinson's constraint qualification  holds at $\bar x$   if  there exists $0\neq \nu \in \R^n$ such that  for each $j \in I(\bar x)$  and  each $\eta_j \in \p_T g_j(\bar x),$    one has $\lan \eta_j, \nu \ran < 0,$  where $\p_T g_j(\bar x)$ is the tangential subdifferential of $g_j$ at $\bar x.$
\end{defi}

\mes
\begin{defi} \label{def2.3} {\bf $($Non-smooth Version of   Abadie's Constraint Qualification $(NACQ)).$} Let  $K =\{x \in \R^n : g_j(x) \le 0, \ \forall \ j=1,2,\cdots, m\}$ be  as in $($\ref{4.1}$),$  and let $C$ be a non-empty closed convex subset of $\R^n$ such that $C \cap K \neq \emptyset.$ Let $\tK := C \cap K,$  and   $\bar x\in \tK.$    We say that    non-smooth  Abadie's constraint qualification  holds at $\bar x$   if  $D(\bar x) \subseteq T_{\tK}(\bar x).$  \end{defi}

\mes

\noin Obviously, the  above definitions of  non-smooth version of constraint qualifications  reduce to their counterparts in  the case of differentiability \cite{ba, bor}.

\mes
\noin  Clearly,  in view of Lemma \ref{lem10}, Definition \ref{def2.2}  and  Definition \ref{def2.3},  the following  implication holds.
\be \label{6666} (NRCQ) \Longrightarrow (NACQ). \ee

\mes
\noin The following example shows that  non-smooth Abadie's constraint qualification is  weaker than non-smooth  Robinson's Constraint Qualification.

\mes
\noin \begin{ex} \label{ex2.1} Let $g_1, g_2, g_3 : \R^2 \lrar \R$ be defined by \beqa &&g_1(x_1, x_2):= |x_2| - x_1, \ g_2(x_1, x_2) := 1-x_1^2 -(x_2 - 1)^2, \ \mbox{and} \\ && g_3(x_1, x_2) :=  1-x_1^2 -(x_2 + 1)^2, \  \forall \ (x_1, x_2) \in \R^2. \eeqa
Then, we have \beqa K &=& \{(x_1, x_2) \in \R^2 : g_j(x_1, x_2) \le 0, \ \forall \ j=1,2,3 \} \\ &=&\{(x_1, x_2) \in \R^2 : 0 \le x_2 \le x_1 \} \cup \{(x_1, x_2) \in \R^2 : 0 \le - x_2 \le x_1 \}. \eeqa
\noin Let \beqa C:=\{(x_1, x_2) \in \R^2 : x_1 \ge 0 \}, \eeqa and $\tK :=C \cap K = K.$  Let $\bar x:=(0, 0) \in \tK.$  Clearly, $g_1, g_2, g_3$ are tangentially convex at $\bar x,$  and $g_1(\bar x) = g_2(\bar x) = g_3(\bar x) = 0.$   Moreover, it is easy to check that \beqa &&g_1^\prime(\bar x, (t_1, t_2)) =|t_2| - t_1, \ g_2^\prime(\bar x, (t_1, t_2)) = 2 t_2, \ \mbox{and} \\ && g_3^\prime(\bar x, (t_1, t_2)) = - 2 t_2, \ \forall \ (t_1, t_2) \in \R^2. \eeqa
This together with $(\ref{42})$  implies that \beqa &&\p_T g_1(\bar x) = \co\{(-1, -1), (-1, 1) \}, \ \p_T g_2(\bar x) =\{(0, 2)\}, \ \mbox{and} \\ &&
\p_T g_3(\bar x) = \{(0, -2)\}. \eeqa  It is clear that  non-smooth Robinson's constraint qualification does not hold at $\bar x.$ But, we have \beqa  T_{\tK}(\bar x) = D(\bar x) = \{(t_1, 0) \in \R^2 : t_1 \ge 0 \}, \eeqa  and hence,  non-smooth Abadie's constraint qualification holds at $\bar x.$
\end{ex}

\mes
\begin{rem}  \label{rem123}  It is worth noting that in \cite{ms}  it has been given a characterization of constrained best approximation under smooth Robinson's constraint qualification  with   continuously Fr\'{e}chet differentiable  constraints, while in this paper we give characterizations of constrained best approximation under non-smooth Abadie's constraint qualification with  the constraint functions are only tangentially convex at  the point of best approximation. Since, in view of $($\ref{6666}$)$  and Example \ref{ex2.1},  $(NACQ)$ is weaker than $(NRCQ).$  So,  our results  are stronger and extend the  obtained results in \cite{ms} and the corresponding results  of \cite{DE,DLW1,DLW2,jm1,jm2,LJ,LN1,LN2}.
\end{rem}

\mes

\noin  For a non-empty subset $W$ of $\R^n$ and an arbitrary point $x\in \R^n,$ we define
\beqa d(x,W):=\inf_{w\in W}\Vert x-w\Vert . \eeqa We say that a point $x_0 \in W$ is a best approximation $($a projection$)$ of $x \in \R^n$ if $\|x - x_0\| = d(x, W)$ \cite{S2}.   The set of all best approximations $($projections$)$ of $x$ in $W$ denoted by $P_W(x)$ and is given by: \beqa P_W(x) :=\{w \in W : \|x - w\| = d(x, W) \}. \eeqa

\mes
\noin The following characterization of best approximation in Hilbert spaces is well known \cite{DB}.

\mes

\begin{lem} \label{lem6.0}
Let $D$ be a non-empty closed convex subset of a Hilbert space $H,$  $x \in H$ and $x_0 \in D.$ Then, $
x_0 = P_D(x)$ if and only if $x - x_0 \in (D - x_0)^\circ.$
\end{lem}

\mes
\noin In the following, we give the notion of strong CHIP.  The definition of strong CHIP was first introduced in \cite{DLW2} (see also, \cite{DE,DB}).

\mes
\begin{defi} \label{d10} {\bf $($Strong CHIP}$).$  Let $C_1, C_2,\ldots, C_m$  be  non-empty  closed convex sets in $\R^n,$ and let $x \in  \bigcap_{j=1}^{m} C_j.$ Then, the collection $\{C_1, C_2,\ldots, C_m\}$  is said to have the strong CHIP $($canonical hull intersection property$)$ at $x$ if \beqa \bigg(\bigcap_{j=1}^{m} C_j  - x\bigg)^\circ = \sum_{j=1}^{m} \bigg(C_j - x\bigg)^\circ. \eeqa
\noin The collection $\{C_1, C_2,\ldots, C_m\}$  is said to have the strong CHIP if it has the strong CHIP at each $x \in \bigcap_{j=1}^{m} C_j.$ \end{defi}

\mes
\noin We recall \cite{BS,bor} the following well known result from the non-smooth analysis.

\mes
\begin{thm} \label{thm00} Let $C \subset \R^n$  be a non-empty  convex set,  and  let   $f : \R^n \lrar \R \cup \{+\infty\}$  be  a proper convex  function such that ${C \cap \dom(f)} \neq \emptyset.$  Assume that   $\bar x \in  C$ and $f$   is  continuous  at $\bar x.$     Then, $\bar x$  is a global  minimizer of the function  $f$ over $C$ if and only if \beqa 0 \in \p f(\bar x) + N_C(\bar x), \eeqa  where domain of the function $f,$ $\dom(f),$ is defined by \beqa \dom(f) :=\{x \in \R^n : f(x) < +\infty \}. \eeqa \end{thm}

\mes

\mes

\section{Dual Cone Characterizations of the Constraint Set $K$}
\noin In this section, we give  dual cone characterizations of  the constraint  set $K$ at a point $x\in K,$ where $K$ is given by (\ref{4.1}). Also, we present sufficient conditions for which the strong CHIP property holds.

\mes

\noin For each $x \in K,$ put \be \label{4.2} M(x) : = \bigcup\limits_{\gl \in S} \bigg\{\sum_{j=1}^{m} \gl_j \p_T g_j(x): \gl_j g_j(x) = 0, \ j=1,2,\cdots,m \bigg\}, \ee where $\gl :=(\gl_1,\gl_2,\cdots,\gl_m) \in S$ and $S :=\R_+^m.$

\mes
\begin{rem} Throughout the paper, we assume that the constraint functions $g_j,$  $j=1,2,\ldots,m,$  are tangentially convex at a given point $\bar x \in \tK :=C \cap K.$  \end{rem}

\mes
\noin  We now give  a dual cone characterization of the constraint set $K,$ which has a crucial role for   characterizing best approximations by the set $\tK := C \cap K.$  Note that $K$ is   not necessarily a closed or a convex set.

\mes

\begin{thm} \label{thm1}  Let $K$ be closed, given by $(\ref{4.1}),$    and let $C$ be a non-empty closed convex  subset of $\R^n$ such that $C \cap K \neq \emptyset.$  Let $\tK := C \cap K,$  $\bar x \in \tK$ and  $M(\bar x)$ be as in $(\ref{4.2}).$    Assume that $K$  is nearly convex at the point $\bar x.$  If   non-smooth  Abadie's  constraint qualification holds at  $\bar x,$  then, $M(\bar x) = (K - \bar x)^\circ = (\tK - \bar x)^\circ.$    \end{thm}
\pr  It is easy to see that $(K - \bar x)^\circ \subseteq (\tK - \bar x)^\circ.$   Now, let $u \in (\tK- \bar x)^\circ$ be arbitrary.  Then, $\lan u, y- \bar x \ran \le 0$ for all $y \in \tK,$
 and so,  \be \label{q000} \lan -u, y - \bar x \ran \ge 0, \  \forall \ y \in \tK. \ee Let $h:\R^n \lrar \R$ be  defined by
 \be \label{102} h(y) := \lan -u, y \ran, \ \forall \ y \in \R^n.  \ee   It is clear that $h$ is a continuous convex function on $\R^n.$
Now, we show that $h(y) \ge 0$ for all $y \in T_{\tK}(\bar x).$  To this end, let $y \in T_{\tK}(\bar x)$ be arbitrary. Then by (\ref{1.2})  there exist $\{t_m\}_{m \ge 1} \subset \R_{++}$ and $\{y_m\}_{m \ge 1}  \subset \R^n$ such that $t_m \lrar 0^+,$  $y_m \lrar y$  and  $\bar x + t_m y_m \in \tK$  for all $m \ge 1.$  Thus, in view of (\ref{q000}), we conclude that   $\lan -u, y_m \ran \ge 0$ for all $m \ge 1.$  This together with  the fact that $y_m \lrar y$ implies that $\lan -u, y \ran \ge 0,$ and so,  \be \label{q100} h(y) \ge 0, \  \forall \ y \in T_{\tK}(\bar x). \ee \noin Consider the following optimization problem: \be \label{p} \min h(y) \ \mbox{subject to} \ y \in T_{\tK}(\bar x).  \ee
It follows from (\ref{q100}) that $y=0 \in T_{\tK}(\bar x)$ is a  global  minimizer of the problem (\ref{p}) over $T_{\tK}(\bar x).$  On the other hand, since $K$  is nearly convex at  the point $\bar x,$  in view of Lemma \ref{lem10} and the validity  of   non-smooth Abadie's  constraint qualification  at $\bar x,$  we have, $D(\bar x) = T_{\tK}(\bar x).$  Therefore, the problem (\ref{p}) can be represented as the following convex optimization problem: \be \label{p110}  \min h(y) \ \mbox{subject to} \ y \in D(\bar x). \ee  Note that $D(\bar x)$ is a closed convex subset of $\R^n,$ and $y=0 \in D(\bar x)$ is a   global minimizer of the problem (\ref{p110}) over $D(\bar x).$   In view of (\ref{4.35}) and  (\ref{1.3}), the problem (\ref{p110}) can be represented as the following convex optimization problem: \be \label{p1}  \min h(y) \ \mbox{subject to} \ y \in \R^n, \ \mbox{and} \ \ g_j^\prime(\bar x, y) \le 0, \ \forall \ j \in I(\bar x). \ee
Note that $y=0$ $($because $g_j^\prime(\bar x, 0) = 0, \ j \in I(\bar x))$  is a   global minimizer of the problem (\ref{p1}).  Let \beqa C_j :=\{y \in \R^n : g_j^\prime(\bar x, y) \le 0 \}, \ (j \in I(\bar x)), \eeqa and let $H := \cap_{j \in I(\bar x)} C_j.$  Since $g_j^\prime(\bar x, \cdot)$ is convex on $\R^n$ $(j \in I),$  it is easy to see that $C_j$ is convex for each $j \in I(\bar x),$ and hence, $H$ is a convex set.  Since, by (\ref{p110}),   $y=0 \in H$ $($note that $g_j^\prime(\bar x, 0) =0 \le 0, \ j \in I(\bar x))$  is a  global  minimizer of the problem (\ref{p1}) over $H,$  it follows from Theorem \ref{thm00} that \be \label{q1} 0 \in \p h(0) + N_H(0). \ee This together with \cite[Section 3.3]{bor}  implies that  \be \label{q2} 0 \in \p h(0) + \sum_{j \in I(\bar x)} N_{C_j}(0). \ee \noin Let \be \label{p3} M_j(\bar x) := \{\gl_j \eta_j : \gl_j \ge 0, \ \eta_j \in \p g_j^\prime(\bar x, \cdot)(0) \}, \ (j \in I(\bar x)). \ee  It is easy to check that $M_j(\bar x)$ is a closed convex cone in $\R^n$ for each $j \in I(\bar x).$  Now, we claim that \be \label{qqq1} N_{C_j}(0) \subseteq M_j(\bar x), \ (j \in I(\bar x)). \ee
\noin Assume if possible that there exists $x^* \in N_{C_j}(0)$ such that $x^* \notin M_j(\bar x).$  Since $M_j(\bar x)$ is a closed convex cone, by using the separation theorem there exists $0 \neq \nu \in \R^n$ such that \be \label{p4} \lan \nu, u_j \ran \le 0 < \lan \nu, x^* \ran, \ \forall \ u_j \in M_j(\bar x), \ (j \in I(\bar x)). \ee

\mes
\noin For simplicity, put $h_j(\cdot) : = g_j^\prime(\bar x, \cdot)$ $(j \in I(\bar x)).$  Since $g_j$ is tangentially convex at $\bar x,$ it follows that, for each $j \in I(\bar x),$    $h_j$ is a  real valued positively homogeneous and  convex function, and so,   \beqa  h_j^\prime(0, \nu) = \max_{\eta_j \in \p h_j(0)} \lan \eta_j, \nu \ran, \ (j\in I(\bar x)). \eeqa  This together with (\ref{p4}) implies that $h_j^\prime(0, \nu) \le 0$ $(j \in I(\bar x)).$   On the other hand, since $h_j$ $(j \in I(\bar x))$ is positively homogeneous, we conclude that $h_j^\prime(0, \nu) = g_j^\prime(\bar x, \nu)$ $(j \in I(\bar x)).$  So, $g_j^\prime(\bar x, \nu) \le 0$ $(j \in I(\bar x)),$ and hence, $\nu \in C_j$ $(j \in I(\bar x)).$   But, we have $x^* \in N_{C_j}(0)$  $(j \in I(\bar x)).$ Therefore, $\lan \nu, x^* \ran \le 0,$ which contradicts (\ref{p4}).  Then, the inclusion (\ref{qqq1}) holds.  This together with (\ref{q2}) implies that \beqa 0 \in \p h(0) + \sum_{j \in I(\bar x)} M_j(\bar x).   \eeqa  So, for each $j \in I(\bar x),$    there exists $\gl_j \ge 0$  such that \be \label{p6} 0 \in \p h(0) + \sum_{j \in I(\bar x)}  \gl_j \p g_j^\prime(\bar x, \cdot)(0).  \ee  It is not difficult to see that $\p g_j^\prime(\bar x, \cdot)(0) = \p_T g_j(\bar x)$ $(j \in I(\bar x)).$  Thus, it follows from (\ref{p6}) that \be \label{p7} 0 \in \p h(0) + \sum_{j \in I(\bar x)} \gl_j \p_T g_j(\bar x). \ee

\mes
\noin But, in view of (\ref{102}), we have $\p h(0) = \{-u\}.$   Now, for each $j \notin I(\bar x),$ put $\gl_j = 0.$  Therefore, we obtain from (\ref{p7}) that \beqa  u \in  \sum_{j \in I(\bar x)} \gl_j \p_T g_j(\bar x) = \sum_{j=1}^{m} \gl_j \p_T g_j(\bar x), \ \mbox{and} \ \ \gl_j g_j(\bar x) = 0, \ j=1,2,\cdots,m, \eeqa  and so, by (\ref{4.2}),  $u \in M(\bar x).$  Hence, $(\tK- \bar x)^\circ \subseteq M(\bar x).$

\mes
\noin Now, we show that $M(\bar x) \subseteq (K - \bar x)^\circ.$ To this end,   let $u \in M(\bar x)$ be arbitrary.
Then, in view of (\ref{4.2}),  there exists $(\gl_1, \gl_2, \cdots, \gl_m) \in S$  with $\gl_j g_j(\bar x) = 0$ $(j =1,2,\cdots,m)$  such that \beqa \label{4.3} u \in   \sum_{j=1}^{m} \gl_j \p_T g_j(\bar x).\eeqa  This implies that, for each $j =1,2,\cdots,m,$   there exists $\eta_j \in \p_T g_j(\bar x)$ such that \be \label{p9} u = \sum_{j=1}^{m} \gl_j \eta_j. \ee
Now, let $y \in K$ be arbitrary. Since   $\bar x \in K$ and $K$  is nearly convex  at  $\bar x,$  it follows  from Definition  \ref{def00} that  there exists a sequence $\{\ga_k\}_{k \ge 1} \subset \R_{++}$ with $\ga_k \lrar 0^{+}$  such that  $\bar x + \ga_k (y- \bar x) \in K$  for all  sufficiently large $k \in \mathbb{N}.$  So, by (\ref{4.1}), \be \label{1000} g_j(\bar x+ \ga_k (y- \bar x)) \le 0, \  \mbox{for all sufficiently large} \ k \in \mathbb{N} \ \mbox{and all} \ j=1,2,\cdots,m.  \ee
 Since $g_j$ $(j=1,2,\cdots,m)$  is tangentially convex at $\bar x,$  it follows from (\ref{4.35}), (\ref{p9}), (\ref{1000}) and the fact that $\gl_j = 0$ for each $j \notin I(\bar x)$ (because $\gl_j g_j(\bar x) =0, \ j=1,2,\ldots,m)$  that
\beqa \lan u, y- \bar x\ran & =& \lan \sum_{j=1}^{m} \gl_j  \eta_j, y- \bar x \ran \\ &=& \sum_{j=1}^{m} \gl_j  \lan  \eta_j, y- \bar x \ran \\ &\le& \sum_{j =1}^{m} \gl_j g_j^\prime(\bar x, y - \bar x) \nonumber \\ &=& \sum_{j \in I(\bar x)} \gl_j g_j^\prime(\bar x, y - \bar x)
  \nonumber \\ &=& \sum_{j \in I(\bar x)} \gl_j \bigg\{
\lim_{k \lrar +\infty} \frac{g_j(\bar x+ \ga_k (y- \bar x)) -  g_j(\bar x)}{\ga_k}\bigg\} \nonumber  \\ &=& \sum_{j \in I(\bar x)} \gl_j \bigg\{\lim_{k \lrar +\infty} \frac{g_j(\bar x+ \ga_k (y- \bar x))}{\ga_k}\bigg\} \nonumber \\ &\le& 0, \ \forall \ y \in K.\eeqa Hence, $u \in (K - \bar x)^\circ,$  which completes the proof.  \hfill \rule{2mm}{2mm}

\mes

\begin{rem} It should be noted that in Theorem \ref{thm1}, for proving the inclusion $M(\bar x) \subseteq (K - \bar x)^\circ,$ we  only require that the set $K$  is nearly convex at  the point $\bar x$  without the validity  of  non-smooth  Abadie's  constraint qualification at  $\bar x.$   \end{rem}

\mes

\begin{cor} \label{p5.6}  Let $K$ be  closed, given by  $(\ref{4.1}),$  and let  $C$  be a non-empty closed convex subset of $\R^n$  such that $C \cap K \neq \emptyset.$  Let $\tK:= C \cap K,$  and $\bar x \in \tK.$   Assume that $K$  is nearly convex at the point $\bar x$  and    non-smooth Abadie's  constraint qualification holds at  $\bar x.$  Then, $\{C, K\}$ has  the strong CHIP at $\bar x.$ \end{cor}
\pr We first note that one always has \beqa (C - \bar x)^\circ + (K - \bar x)^\circ \subseteq (\tK - \bar x)^\circ. \eeqa \noin For the converse inclusion, in view of  Theorem \ref{thm1}, we conclude that  $(\tK - \bar x)^\circ = (K - \bar x)^\circ.$  Therefore, since $0 \in (C - \bar x)^\circ,$  we have \beqa (\tK - \bar x)^\circ \subseteq (C - \bar x)^\circ + (K - \bar x)^\circ, \eeqa  which completes the proof.   \hfill  \rule{2mm}{2mm}

\mes

\noin The following example shows that   the converse statement to Corollary \ref{p5.6} is not valid.

\mes
\begin{ex} \label{ex3.1} Let $g(x):=|x| - x$ for all $x \in \R.$   Thus, $K=\{x\in\R  :  g(x)\le 0\}= [0, +\infty),$ which is closed.    Let $C:= (-\infty, 0],$ and  $\bar x = 0.$ It is clear that $g$ is tangentially convex at $\bar x,$ $g^\prime(\bar x, \nu) = |\nu| - \nu$ for all $\nu \in \R,$  and $C \cap K = \{0\}.$  Also, we have  $g(\bar x) =0$ and $\p_T g(\bar x) = [-2, 0].$ Note that $K$  is nearly convex at  the point $\bar x.$   It is easy to check that $(K - \bar x)^\circ = (-\infty, 0] $ and $(C - \bar x)^\circ = [0, +\infty).$   Let $\tK : = C \cap K = \{0\}.$  Then, \beqa (\tK - \bar x)^\circ  &&= \{0\}^\circ \\ &&= \R  \\ &&= (-\infty, 0] + [0, +\infty) \\ &&= (K - \bar x)^\circ + (C - \bar x)^\circ. \eeqa Thus, $\{C, K\}$ has the strong CHIP at $\bar x.$ But, on the other hand, we have \beqa  D(\bar x) = [0, +\infty) \ \mbox{and} \  \ T_{\tK}(\bar x) = \{0\}.  \eeqa  This implies that \beqa  D(\bar x) \nsubseteq T_{\tK}(\bar x), \eeqa and hence,   non-smooth Abadie's  constraint qualification  does not hold at  $\bar x.$
\end{ex}

\mes
\noin By the following example we show   that  in Theorem \ref{thm1} the validity of  the nearly convexity of $K$   at the point $\bar x \in K$  cannot be omitted.

\mes
\noin \begin{ex} Let $g_1, g_2: \R \lrar \R$ be defined by  \beqa g_1(x) := 8-x^3, \ \mbox{and} \ \ g_2(x) := -x^2 + 6 x - 8, \ \forall \ x \in \R. \eeqa  Then,  we have  \beqa K :=\{x \in \R : g_j(x) \le 0, \ j=1,2 \} = \{2\} \cup [4, +\infty), \eeqa  which is closed.   Let $C:=[2, 3],$   $\tK := C \cap K = \{2\}$ and $\bar x :=2.$  It is easy to see that $g_1$ and $g_2$ are tangentially  convex at $\bar x,$    $g_1(\bar x) = g_2(\bar x) =0,$   $\p_T g_1(\bar x) =\{-12\}$ and $\p_T g_2(\bar x) = \{2\}.$   Clearly, \beqa D(\bar x) = \{0\}, \ \mbox{and} \ \ T_{\tK}(\bar x) = \{0\}. \eeqa  Thus, $D(\bar x) = T_{\tK}(\bar x),$ and so,   non-smooth  Abadie's  constraint qualification  holds at  $\bar x,$  while it is clear that $K$  is not nearly convex at  the point $\bar x.$

\mes

\noin On the other hand, it is not difficult to check that $M(\bar x) = \R,$  $(\tK - \bar x)^\circ = \R$  and  $(K - \bar x)^\circ = (-\infty, 0].$   Hence,  Theorem \ref{thm1} does not hold. \end{ex}

\mes

\noin The  following example shows that   non-smooth Abadie's  constraint qualification  in Theorem \ref{thm1} cannot be omitted.

\mes

\begin{ex}  \label{ex3.4}  Let $g_1, g_2: \R \lrar \R$ be defined by
\beqa g_1(x):=\left \{ \begin{array}{ll}  {x}^{\frac{3}{2}}, & x>0,
 \\0, & x \le 0, \end{array}\right. \eeqa  and
\beqa g_2(x):=\left \{ \begin{array}{ll} - {x}^{\frac{3}{2}}, & x>0,
 \\0, & x \le 0. \end{array}\right. \eeqa  Thus, we have \beqa K :=\{x \in \R : g_j(x) \le 0, \ j=1,2 \} = (-\infty, 0], \eeqa which is closed.   Let $C:=[0, 1],$  $\tK :=C \cap K = \{0\}$ and $\bar x := 0.$    It is easy to check that $g_1$ and $g_2$ are tangentially  convex at $\bar x,$  $g_1(\bar x) = g_2(\bar x) =0,$  $\p_T g_1(\bar x) = \{0\}$ and $\p_T g_2(\bar x) = \{0\}.$  Also, one can see that  \beqa D(\bar x) = \R, \ \mbox{and} \ \ T_{\tK}(\bar x) = \{0\}. \eeqa  Therefore, $D(\bar x) \nsubseteq T_{\tK}(\bar x),$ and hence,   non-smooth Abadie's  constraint qualification  does not hold  at  $\bar x,$  while $K$   is nearly convex  at  the point $\bar x.$

\mes

\noin Furthermore, it is easy to see that ${M}(\bar x) =\{0\},$  $(\tK - \bar x)^\circ = \R$  and $(K - \bar x)^\circ = [0, +\infty),$ and so,  Theorem  \ref{thm1} does not hold.
\end{ex}

\mes

\mes

\section{Characterizations of Constrained Best Approximation}
\noin In this section,   we give characterizations of constrained best approximations under   non-smooth Abadie's  constraint qualification. Let $K$ be as in (\ref{4.1}), given by,   \beqa  \label{6.1} K:= \{x \in \R^n : g_j(x) \le 0, \ j=1,2,\cdots,m \}, \eeqa  where $g_j : \R^n \lrar \R$ $(j=1,2, \cdots,m)$  is a tangentially convex function at a given  point $\bar x \in K.$   Let  $S := \R_+^m,$  and let  $C$  be a non-empty closed convex  subset of $\R^n$ such that $C \cap K \neq \emptyset.$ Note that  $K$ is  not necessarily a closed or a convex set.

\mes
\noin The following theorem shows that under  non-smooth Abadie's  constraint qualification  the "perturbation property" is characterized by the strong conical hull intersection property (Strong CHIP).

\mes

\begin{thm} \label{thm6.1}  Let $K$ be  closed, given by $(\ref{4.1}),$  and  let $C$  be a non-empty closed convex  subset of $\R^n$ such that $C \cap K \neq \emptyset.$  Let  $\bar x \in \tK:=C\cap K.$   Assume that  $\tK$ is  closed and convex.  If $K$  is nearly convex at  the point $\bar x$  and   non-smooth Abadie's  constraint qualification holds at  $\bar x,$  then the following assertions are equivalent. \\
\noin $(i)$ $\{C, K\}$ has the strong CHIP at $\bar x,$ \\
\noin $(ii)$ For any $x \in \R^n,$ $\bar x = P_{\tK}(x)$ if and only if there exist  $(\gl_1, \gl_2, \cdots, \gl_m) \in S$  with  $\gl_j g_j(\bar x) =0,$ $j=1,2,\cdots,m,$ and  \be \label{6.2}
 \eta_j \in \p_T g_j(\bar x) \ (j=1,2,\cdots,m) \  \mbox{such that} \ \   \bar x = P_C(x - \sum_{j=1}^{m} \gl_j \eta_j).  \ee \end{thm}
\pr  $[(i)\Longrightarrow (ii)].$  Suppose that $(i)$ holds. Then, by  Definition \ref{d10}, \be \label{6.3} (\tK-\bar x)^\circ = (C-\bar x)^\circ + (K-\bar x)^\circ.\ee  Also, in view of the hypotheses and Theorem \ref{thm1}, we have $M(\bar x) = (K - \bar x)^\circ.$  So, it follows from (\ref{6.3}) that \be \label{6.4} (\tK-\bar x)^\circ = (C-\bar x)^\circ + M(\bar x).\ee  \noin Now,  for any $x \in \R^n,$  assume that $\bar x = P_{\tK}(x).$  Thus, by
Lemma \ref{lem6.0}, one has $x - \bar x \in (\tK - \bar x)^\circ.$  Therefore, in view of  (\ref{4.2}) and (\ref{6.4}), there exist $\ell \in (C-\bar x)^\circ$ and $(\gl_1, \gl_2, \cdots, \gl_m) \in S$  with $\gl_j g_j(\bar x) =0$ $(j=1,2,\cdots,m)$  such that
\beqa x-\bar x - \ell \in \sum_{j=1}^{m} \gl_j \p_T g_j(\bar x). \eeqa  So, for each $j=1,2,\cdots,m,$ there exists $\eta_j \in \p_T g_j(\bar x)$ such that  \beqa x-\bar x - \ell = \sum_{j=1}^{m} \gl_j \eta_j. \eeqa   Then, we conclude  that \begin{eqnarray*}
[x - \sum_{j=1}^{m} \gl_j \eta_j] - \bar x= \ell \in (C-\bar x)^\circ,  \end{eqnarray*} for some
$(\gl_1, \gl_2, \cdots, \gl_m) \in S$   with  $\gl_j  g_j(\bar x) = 0,$ and some  $\eta_j \in \p_T g_j(\bar x), \ j=1,2,\cdots,m.$  Frow now on, by an argument similar  to the proof of Theorem 4.1 $($the implication $[(i) \Longrightarrow (ii)])$  in \cite{jm3} and using Theorem \ref{thm1} the result follows. \\
\noin $[(ii)\Longrightarrow (i)].$ Let $y \in (\tK-\bar x)^\circ$ be arbitrary, and let $x:= \bar x + y.$ Thus, by Lemma \ref{lem6.0},
$\bar x  =P_{\tK}(x).$  Then, in view of the hypothesis $(ii),$ there exist  $(\gl_1, \gl_2, \cdots, \gl_m) \in S$  with  $\gl_j g_j(\bar x) =0,$ $j=1,2,\cdots,m,$ and  \beqa
 \eta_j \in \p_T g_j(\bar x) \ (j=1,2,\cdots,m)
 \   \mbox{such that}  \ \
 \bar x = P_C(x - \sum_{j=1}^{m} \gl_j \eta_j). \eeqa  Again, by  Lemma \ref{lem6.0}, \be \label{4b}
 y - \sum_{j=1}^{m} \gl_j \eta_j \in (C-\bar x)^\circ,  \ee  for some
$(\gl_1, \gl_2, \cdots, \gl_m) \in S$   with  $\gl_j  g_j(\bar x) = 0,$ and some $\eta_j \in \p_T g_j(\bar x), \ j=1,2,\cdots,m.$  Now, by an argument similar to the proof of Theorem 4.1 $($the implication $[(ii) \Longrightarrow (i)])$  in \cite{jm3}  and using Theorem \ref{thm1} the proof is completed.
\hfill  \rule{2mm}{2mm}

\mes

\noin The following examples illustrate Theorem \ref{thm6.1}. Moreover, these examples  justify how  one can use best approximations to check the strong CHIP without explicitly proving the strong CHIP property.

\mes

\begin{ex} \label{ex4.1} Let $g_1, g_2:\R^2 \lrar \R$ be defined by \beqa &&g_1(x_1, x_2) := |x_2|-x_1-x_1^2-x_1^3, \ \mbox{and} \\  && g_2(x_1, x_2):= |x_1 - x_2| -x_1-x_1 x_2 - x_2^3,  \ \forall \ (x_1, x_2) \in \R^2. \eeqa  Let $S:= \R_+^2$ and  $C:= \R^2.$   It is easy to see that \beqa K :=\{(x_1, x_2) \in \R^2 : g_j(x_1, x_2) \le 0, \ j=1,2 \} = \{(x_1, x_2) \in \R^2 : x_1 \ge x_2 \ge 0 \}, \eeqa  which is closed and convex.   Let $\tK := C \cap K = K$ and $\bar x :=(0, 0) \in \tK.$  Note that $\tK$ is closed and convex. It is clear that $g_1$ and $g_2$  are tangentially convex at $\bar x,$  but  not convex. Moreover,  $g_1(\bar x) = g_2(\bar x) = 0,$ and \beqa g_1^\prime(\bar x, (t_1, t_2)) = |t_2|-t_1, \ \mbox{and} \ \ g_2^\prime(\bar x, (t_1, t_2)) = |t_1 -t_2| - t_1, \ \forall \ (t_1, t_2) \in \R^2. \eeqa  Therefore, we have \beqa \p_T g_1(\bar x) = \co \{(-1, 1), (-1, -1)\}, \ \mbox{and} \ \ \p_T g_2(\bar x) = \co \{(-2, 1), (0, -1)\}. \eeqa So, it is easy to check that  \beqa D(\bar x) =\{(t_1, t_2) \in \R^2 : t_1 \ge t_2 \ge 0\} = T_{\tK}(\bar x), \eeqa  and hence,   Abadie's  constraint qualification holds at  $\bar x.$ Also, it is clear that $K$  is nearly convex  at the point $\bar x.$

\mes

\noin Now, for any $x:=(0, x_2) \in \R^2$ with $x_2 \le 0,$  it is easy to see that \beqa P_{\tK}(x) = \bar x = (0,0) =P_C((0,0)) = P_C(x - (\gl_1 \eta_1 + \gl_2 \eta_2)), \eeqa where $(\gl_1 :=0,  \gl_2 := - x_2) \in S,$    $\gl_j g_j(\bar x)  = 0$ $(j=1,2),$  and $\eta_1:=(-1, 1) \in \p_T g_1(\bar x),$  $\eta_2 :=(0, -1) \in \p_T g_2(\bar x).$   Then, in view of Theorem \ref{thm6.1} $($the implication $[(ii)\Longrightarrow (i)]),$  we conclude that  $\{C, K\}$ has the strong CHIP at $\bar x.$  Indeed, one can see \beqa (\tK - \bar x)^\circ &=& \{(t_1, t_2) \in \R^2 : t_1 \le -t_2, \ t_1 \le 0 \} \cup (\R_{-} \times \R_{-}) \\ &=& \{(0, 0)\} + \{(t_1, t_2) \in \R^2 : t_1 \le -t_2, \ t_1 \le 0 \} \cup (\R_{-} \times \R_{-}) \\ &=& (C - \bar x)^\circ + (K - \bar x)^\circ. \eeqa \end{ex}

\mes

\begin{ex}  \label{ex4.2} Let $g_1, g_2 : \R \lrar \R$ be defined by \beqa g_1(x) := 1-x^3, \ \mbox{and} \ \ g_2(x) := x^3 -3x^2 +  x -3, \ \forall \ x \in \R. \eeqa Let $S:= \R_+^2$ and $C:=[1, +\infty).$    Clearly, we have \beqa K =\{x \in \R : g_j(x) \le 0, \ j=1,2\} = [1, 3],\eeqa  which is closed and convex.    Let $\tK :=C \cap K = [1, 3]$ and $\bar x := 1 \in \tK.$  Thus, $\tK$ is closed and convex, $g_1, g_2$ are tangentially convex at $\bar x$ $($but  not convex$),$
$g_1(\bar x) = 0,$ $g_2(\bar x) = -4 \neq 0,$  and  \beqa g_1^\prime(\bar x, t) = -3t, \ \mbox{and} \ \ g_2^\prime(\bar x, t) = -2t, \ \forall \ t \in \R. \eeqa  This implies that \beqa \p_T g_1(\bar x) =\{-3\}, \ \mbox{and} \ \ \p_T g_2(\bar x) = \{-2\}. \eeqa
Moreover, it is not difficult to show that   \beqa D(\bar x) = [0, +\infty) = T_{\tK}(\bar x), \eeqa and so,  non-smooth Abadie's  constraint qualification holds at  $\bar x.$ Note that $K$  is nearly convex at the point $\bar x.$   It is easy to see that, for any $x \in \R$ with $x \le 1,$  we have \beqa P_{\tK}(x) = \bar x = 1 = P_C(1) = P_C(x - (\gl_1 \eta_1 + \gl_2 \eta_2)), \eeqa where $(\gl_1 := \frac{1-x}{3}, \gl_2 :=0) \in S,$  $\gl_j g_j(\bar x) =0$ $(j=1,2),$  and $\eta_1 := -3 \in \p_T g_1(\bar x),$ $\eta_2 :=-2 \in \p_T g_2(\bar x).$   Hence, by using Theorem \ref{thm6.1} $($the implication $[(ii)\Longrightarrow (i)]),$  we conclude that  $\{C, K\}$ has the strong CHIP at $\bar x.$   Indeed, one can see \beqa (\tK - \bar x)^\circ &=& (-\infty, 0]  \\ &=& (-\infty, 0] + (-\infty, 0]  \\ &=& (C - \bar x)^\circ + (K - \bar x)^\circ. \eeqa
\end{ex}

\mes

\noin Now, let $x \in \R^n$ be fixed, and define the function $h:\R^n \lrar [0, +\infty)$ by
\beqa h(y) : =\|y - x\|, \ \forall \ y \in \R^n. \eeqa For $\bar x \in \R^n,$ we recall that $\p h(\bar x):=\p \|\cdot- x\|(\bar x)$ is
 given by \be \label{4.44} \p \|\cdot - x\|(\bar x) = \{x^* \in \R^n : \|x^*\| = 1, \ \lan x^*, \bar x -x \ran = \|\bar x - x\| \}. \ee

\mes
\noin  In the following, we give  the Lagrange multipliers characterization of constrained best approximation under  non-smooth Abadie's  constraint qualification.

\mes
\begin{thm}\label{thm6.2}  Let $K$ be closed, given by  $(\ref{4.1}),$  and  let $C$  be a non-empty closed convex  subset of $\R^n$ such that $C \cap K \neq \emptyset.$  Let  $\bar x \in \tK:=C\cap K$ and $x \in \R^n.$    Assume that  $\tK$ is  closed and convex.  If $K$  is nearly convex at the point $\bar x$  and    non-smooth Abadie's  constraint qualification holds at  $\bar x,$  then the following assertions are equivalent. \\
\noin $(i)$ $\bar x = P_{\tK}(x).$\\
\noin $(ii)$ There exist  $(\gl_1, \gl_2, \ldots, \gl_m) \in S$  with  $\gl_j g_j(\bar x) =0,$ $j=1,2,\cdots,m,$ and  \beqa
 \eta_j \in \p_T g_j(\bar x) \ (j=1,2,\ldots,m) \  \mbox{such that} \ \   \bar x = P_C(x - \sum_{j=1}^{m} \gl_j \eta_j).  \eeqa
\noin $(iii)$  There exist  $(\gl_1, \gl_2, \ldots, \gl_m) \in S$  with  $\gl_j g_j(\bar x) =0,$ $j=1,2,\cdots,m,$ and
$\eta_j \in \p_T g_j(\bar x) \ (j=1,2,\ldots,m)$  such that   \beqa 0 \in \p \|\cdot - x\|(\bar x) + (C-\bar x)^\circ + \sum_{j=1}^{m} \gl_j \eta_j, \eeqa  where we denote  $\p f(x_0)$ for the convex subdifferential of a convex function $f: \R^n \lrar \R$ at the point $x_0 \in \R^n.$
\end{thm}
\pr $[(i) \Longleftrightarrow (ii)].$ Since, by the hypothesis,   non-smooth Abadie's constraint qualification holds at $\bar x$ and $K$  is nearly convex at  the point $\bar x,$  it follows from  Corollary  \ref{p5.6} that  $\{C, K\}$ has the strong CHIP at $\bar x.$  Therefore,  the implication $[(i)\Longleftrightarrow(ii)]$  follows from Theorem \ref{thm6.1}.\\
\noin $[(i) \Longrightarrow (iii)].$ We may assume without loss of generality  that $x \neq \bar x.$ Suppose that $(i)$ holds.   Then, we have $\bar x = P_{\tK}(x).$  This together with Lemma \ref{lem6.0}  implies that $ x - \bar x \in (\tK - \bar x)^\circ.$  But, in view of Theorem \ref{thm1}, one has \beqa (\tK - \bar x)^\circ = M(\bar x). \eeqa
Hence, $x - \bar x \in M(\bar x).$ Since $M(\bar x)$ is a cone, we conclude that \beqa \frac{x - \bar x}{\|\bar x - x\|} \in M(\bar x). \eeqa Therefore, it follows from (\ref{4.2}) that there exist  $(\gl_1, \gl_2, \ldots, \gl_m) \in S$  with  $\gl_j g_j(\bar x) =0,$ $j=1,2,\cdots,m,$ and
$\eta_j \in \p_T g_j(\bar x) \ (j=1,2,\ldots,m)$  such that \be \label{6.100} -u :=\frac{x - \bar x}{\|\bar x - x\|} = \sum_{j=1}^{m} \gl_j \eta_j, \ee
\noin  where \beqa u: = \frac{\bar x - x}{\|\bar x - x\|}. \eeqa  Then, $u \in \R^n,$  $\|u\| = 1$ and \beqa \lan u, \bar x - x \ran = \|\bar x - x\|. \eeqa  This together with (\ref{4.44}) implies that \be \label{6.101} u \in \p \|\cdot - x\|(\bar x). \ee
Thus, it follows from (\ref{6.100}) and (\ref{6.101}) that  \be \label{6.102} 0 \in  \p \|\cdot - x\|(\bar x) + \sum_{j=1}^{m} \gl_j \eta_j, \ee for some $(\gl_1, \gl_2, \cdots, \gl_m) \in S$  with  $\gl_j g_j(\bar x) =0,$  and some
$\eta_j \in \p_T g_j(\bar x),$  $j=1,2,\cdots,m.$  On the other hand, since $0 \in (C - \bar x)^\circ,$  we have \beqa \p \|\cdot - x\|(\bar x) + \sum_{j=1}^{m} \gl_j \eta_j \subseteq \p \|\cdot - x\|(\bar x) + (C - \bar x)^\circ +  \sum_{j=1}^{m} \gl_j \eta_j.  \eeqa Hence, in view of (\ref{6.102}) there exist $(\gl_1, \gl_2, \ldots, \gl_m) \in S$  with  $\gl_j g_j(\bar x) =0,$  and
$\eta_j \in \p_T g_j(\bar x),$ $j=1,2,\cdots,m,$  such that \beqa  0 \in \p \|\cdot - x\|(\bar x) + (C - \bar x)^\circ +  \sum_{j=1}^{m} \gl_j \eta_j,  \eeqa  which implies  $(iii)$ holds. \\
\noin $[(iii)\Longrightarrow (i)].$ Suppose that $(iii)$ holds. Then  there exist  $(\gl_1, \gl_2, \ldots, \gl_m) \in S$  with  $\gl_j g_j(\bar x) =0,$ $j=1,2,\cdots,m,$ and
$\eta_j \in \p_T g_j(\bar x) \ (j=1,2,\ldots,m)$  such that   \be \label{1} 0 \in \p \|\cdot - x\|(\bar x) + (C-\bar x)^\circ + \sum_{j=1}^{m} \gl_j \eta_j. \ee   \noin  Now, let $y \in \tK$ be arbitrary. So, $y \in K.$  Since  $K$  is nearly convex  at  the point $\bar x,$  it follows from  Definition \ref{def00}  that  there exists a sequence $\{\ga_k\}_{k \ge 1} \subset \R_{++}$ with $\ga_k \lrar 0^{+}$  such that  $\bar x + \ga_k (y- \bar x) \in K$  for all  sufficiently large $k \in \mathbb{N}.$  So, by (\ref{4.1}), \be \label{100} g_j(\bar x+ \ga_k (y- \bar x)) \le 0, \  \mbox{for all sufficiently large} \ k \in \mathbb{N} \   \mbox{and all} \ \ j=1,2,\cdots,m.  \ee
 Since $g_j$ $(j=1,2,\cdots,m)$  is tangentially convex at $\bar x,$  it follows from (\ref{4.35}),   (\ref{1})  and (\ref{100}) with some $\nu \in (C - \bar x)^\circ$  that
\begin{eqnarray}  \label{f1}  \|\bar x - x\|
 - \|y - x\|   &\le&  \lan \sum_{j=1}^{m} \gl_j  \eta_j + \nu, y- \bar x \ran \nonumber \\ &=& \la \sum_{j=1}^{m} \gl_j  \eta_j, y- \bar x \ran + \la \nu, y - \bar x \ra \nonumber \\ &\le& \la \sum_{j=1}^{m} \gl_j  \eta_j, y- \bar x \ran  \nonumber \\ &\le& \sum_{j =1}^{m} \gl_j g_j^\prime(\bar x, y - \bar x) \nonumber \\ &=& \sum_{j \in I(\bar x)} \gl_j g_j^\prime(\bar x, y - \bar x)
  \nonumber \\ &=& \sum_{j \in I(\bar x)} \gl_j \bigg\{
\lim_{k \lrar +\infty} \frac{g_j(\bar x+ \ga_k (y- \bar x)) -  g_j(\bar x)}{\ga_k}\bigg\} \nonumber  \\ &=& \sum_{j \in I(\bar x)} \gl_j \bigg\{\lim_{k \lrar +\infty} \frac{g_j(\bar x+ \ga_k (y- \bar x))}{\ga_k}\bigg\} \nonumber \\ &\le& 0, \  \forall \  y \in \tK. \end{eqnarray}
\noin Note that in the above we used the fact that $\gl_j =0$ for each $j \notin I(\bar x),$  because   $\gl_j g_j(\bar x) = 0$  for all  $j=1,2,\ldots,m.$
\noin Therefore,  we conclude from (\ref{f1}) that
 $\|\bar x - x\| = \inf_{y \in \tK} \|y - x\| = d(x, \tK),$ and so, $\bar x = P_{\tK}(x),$ i.e.,  $(i)$ holds. \hfill  \rule{2mm}{2mm}

\mes

\mes

\noin {\bf Acknowledgements.}  This research was partially supported by Mahani Mathematical Research Center, grant no. 97/3267.

\mes

\mes

\end{document}